\documentclass[12pt]{amsart}
\usepackage{mathrsfs, amssymb, array}

\usepackage{mathcomp}
\usepackage[matrix,arrow,curve]{xy}
\usepackage{longtable}

\newcommand{\comment}[1]{}

\newcommand{\OOO}{\mathscr{O}}

\newcommand{\CC}{\mathbb{C}}
\newcommand{\ZZ}{\mathbb{Z}}
\newcommand{\PP}{\mathbb{P}}
\newcommand{\QQ}{\mathbb{Q}}
\newcommand{\FF}{\mathbb{F}}
\newcommand{\RR}{\mathbb{R}}

\newcommand{\Sym}{\mathfrak{S}}

\newcommand{\xref}[1]{{\rm \ref{#1}}}

\newcommand{\comp}{\mathbin{\scriptstyle{\circ}}}
 
\newcommand{\NE}{\overline{\operatorname{NE}}}
\newcommand{\N}{\operatorname{N}}

\newcommand{\Pic}{\operatorname{Pic}}
\newcommand{\Cl}{\operatorname{Cl}}
\newcommand{\Aut}{\operatorname{Aut}}

\newcommand{\rk}{\operatorname{rk}}
\newcommand{\id}{\operatorname{id}}

\newcommand{\GL}{\operatorname{GL}}
\newcommand{\PGL}{\operatorname{PGL}}
\newcommand{\SL}{\operatorname{SL}}

\newcommand{\muu}{{\boldsymbol{\mu}}}

\newcommand{\Fl}{\operatorname{Fl}}
\newcommand{\Nef}{\operatorname{Nef}}
\newcommand{\Exc}{\operatorname{Exc}}

\renewcommand{\xref}[1]{\textup{\ref{#1}}}

\newtheorem{theorem}{Theorem}
\renewcommand{\thetheorem}{\thesection.\arabic{theorem}}
\numberwithin{theorem}{section}
\numberwithin{equation}{theorem}
\renewcommand{\theequation}{\thetheorem.\arabic{equation}}

\newtheorem{mtheorem}[theorem]{}
\newtheorem{stheorem}[equation]{}

\theoremstyle{definition}
\newtheorem{case}[theorem]{}
\newtheorem{scase}[equation]{}

\def\numer{\refstepcounter{equation}{\rm (\theequation)}}%

\title{$G$-Fano threefolds, II}

\author{Yuri Prokhorov}

\thanks{The work was partially supported by the grant RFBR Nos. 11-01-00336-a, 11-01-92613-KO\_a,
Leading Scientific Schools No. 4713.2010.1, 
and 
AG Laboratory HSE RF government grant ag. 11.G34.31.0023}
\address{
Department 
of Algebra, Faculty of Mathematics, Moscow State
University, Moscow, 119 991, RUSSIA
\newline\indent
Laboratory of Algebraic Geometry, SU-HSE, 
7 Vavilova Str., Moscow, 117312, RUSSIA}
\email{prokhoro@gmail.com}
\begin{document}

\begin{abstract}
We classify Fano threefolds with only Gorenstein terminal 
singularities  and Picard number greater than $1$
satisfying  an additional assumption that the 
$G$-invariant part of the Weil divisor class group is of rank $1$
with respect to an action of 
some group $G$.
\end{abstract}

\maketitle

\section{Introduction.}
This work is a sequel to \cite{Prokhorov2010a}.

\begin{case}
\label{notation}
Let $X$ be a Fano threefold with at worst terminal Gorenstein singularities
defined over a field $\Bbbk$ of characteristic $0$. 
Assume that a group $G$ acts on $\bar X:=X\otimes_\Bbbk \bar \Bbbk$, where 
$\bar \Bbbk$ is the algebraic closure of $\Bbbk$.
Moreover, we assume that $X$, $G$ and $\Bbbk$ satisfy
one of the following two conditions.

\begin{itemize}
 \item [(a)] \textit{Algebraic case.} 
$G$ is the Galois group of $\bar \Bbbk$ over $\Bbbk$
acting on $\bar X$ through the second factor.
The action of $G$ on $X$ is trivial.

 \item [(b)] \textit{Geometric case.}
The field $\Bbbk$ is algebraically closed,
$G$ is a finite group, 
and the action of $G$ on $X$ is given by a 
homomorphism $G\to \Aut_{\Bbbk}(X)$.
\end{itemize}
We say that $X$ is a \textit{$G$-Fano threefold} if  $X$ has at worst terminal 
Gorenstein singularities, $-K_X$ is ample, and 
\begin{equation}\label{condition}
\text{$\Cl(X)^G$ is a subgroup of rank $1$
containing $-K_X$.}
\end{equation}
where $\Cl(X)$ is the Weil divisor class group of $X$.
We refer to the introduction in \cite{Prokhorov2009e}, \cite{Prokhorov2010a}, and \cite{Prokhorov2009d}
for the motivation behind this definition.
\end{case}

In this paper we give a classification of one class 
of Gorenstein $G$-Fano threefolds. 
More general we assume that the group $G$ acts only on the Picard lattice 
$\Pic(X)$ in an appropriate way (so we do not assume that $G$ acts on the variety itself).
\begin{mtheorem}{\bf Theorem.}
\label{main}
Let $X$ be a smooth Fano threefold 
over an algebraically closed field of characteristic $0$.
Assume that $\rho(X)>1$ and a finite group $G$ 
acts on $\Pic(X)$ preserving the intersection form and the class $c_1=[-K_X]$.
Furthermore, assume that $\Pic(X)^G\simeq \ZZ$.
Then $X$ is one of 
the varieties in the table below.
\end{mtheorem}

\par\bigskip\bigskip\noindent
\setlongtables
\begin{longtable}
{cccp{232pt}}
\label{ta:res1}
\\
\\
\hline
\\[-5pt]
\multicolumn{1}{c}{No.}
&\multicolumn{1}{c}{$\rho(X)$}
&\multicolumn{1}{c}{$-K_X^3$}
&\multicolumn{1}{c}{$X$}
\\
\\[-10pt]
\hline
\\[-5pt]
\endfirsthead
\\
\hline
\\
\multicolumn{1}{c}{No.}
&\multicolumn{1}{c}{$\rho(X)$}
&\multicolumn{1}{c}{$-K_X^3$}
&\multicolumn{1}{c}{$X$}
\\
\\
\hline
\\
\endhead
\hline
\endlastfoot
\\
\endfoot

\numer\label{list-rho=2-degree=12}
& $2$& $12$ & 
Let $Z_6\subset \PP^8$ be the image of the Segre embedding of $\PP^2\times \PP^2$
and let $Y_6\subset \PP^9$ be the projective cone over $Y_6$. 
Then $X$ is an intersection of $Y_6$ with a hyperplane and a quadric
 \\ [5pt]
\numer\label{list-section-P3P3}& $2$& $20$ & $X$ 
is an intersection of three divisors 
of bidegree $(1,1)$ in $\PP^3\times \PP^3$
 \\ [5pt]
\numer\label{list-2-blowup-Q}& $2$& $28$& 
Let $\sigma: Y'\to \PP^5$ be the blowup with center a Veronese surface.
Then $X$ is an intersection of 
$D_1\in |\sigma^*\OOO_{\PP^5}(1)|$ and $D_2\in |\sigma^*\OOO_{\PP^5}(2)-E|$
where $E$ is the exceptional divisor
 \\ [5pt]
 \numer\label{list-W6}
&$2$& $48$& $V_6\subset \PP^2\times \PP^2$,
a divisor of bidegree $(1,1)$
 \\ [5pt]
\numer\label{list-double-cover-P1P1P1}
& $3$& $12$ & 
 $X$ is a double cover of $\PP^1\times \PP^1\times \PP^1$
 whose branch locus is a member of $|-K_{\PP^1\times \PP^1\times \PP^1}|$
 \\ [5pt]
\numer\label{list-rho=3-degree=30}
& $3$& $30$& 
$X$ is an intersection of divisors 
of tridegrees $(0,1,1)$, $(1,0,1)$,  $(1,1,0)$ in $\PP^2\times \PP^2\times \PP^2$
 \\ [5pt]
\numer\label{list-P1P1P1}& $3$& $48$ & 
 $X=\PP^1\times \PP^1\times \PP^1$
 \\ [5pt]
\numer\label{list-P1P1P1P1} 
& $4$& $24$& $X\subset \PP^1\times \PP^1\times \PP^1\times \PP^1$ is a divisor of 
 multidegree $(1,1,1,1)$
\end{longtable}

\begin{scase}{\bf Remark.}
It is easy to see that all the Fano varieties in the table, except for \eqref{list-rho=2-degree=12} and \eqref{list-double-cover-P1P1P1}, 
are rational. Varieties \eqref{list-rho=2-degree=12} and \eqref{list-double-cover-P1P1P1} are not rational
(see \cite{Alzati1992}).
\end{scase}
\noindent
Theorem \ref{main} will be proved in Sect. \ref{section-Proof-of-Theorem}.
In Sect.
\S \ref{section-Examples-and-remarks} we give several examples.
Sections \ref{section-Action-on-the-Picard-lattice} and 
 \ref{section-Fano-conic-bundles-and-del-Pezzo-fibrations}
are preliminary.
Finally, in  Sect. \ref{section-Singular-case} we investigate singular $G$-Fano threefolds 
with $\rho(X)>1$.

\par\smallskip\noindent
\textbf{Notation.}
We work over an algebraically closed field $\Bbbk$ of 
characteristic $0$.
By $c_1\in \Pic(X)$ we denote the anticanonical class $-K_X$ 
of $X$,   
$\Cl(X)$ denotes the Weil divisor class group.

\par\smallskip\noindent
\textbf{Acknowledgments.}
I would like to thank Sergei Galkin for 
useful comments.

\section{Examples and remarks}\label{section-Examples-and-remarks}
\begin{case} \textbf{Remark.}
By construction, in all cases except for   
\eqref{list-double-cover-P1P1P1} and
\eqref{list-P1P1P1} 
the variety $X$ is a complete 
intersection of certain Cartier divisors 
in some higher-dimensional variety $Y$, where 
$Y=Y_6$ in the case \eqref{list-rho=2-degree=12},
$Y=\PP^3\times \PP^3$ in the case \eqref{list-section-P3P3},
$Y=Y'$ in the case \eqref{list-2-blowup-Q},
$Y=\PP^2\times \PP^2$ in the case \eqref{list-W6},
$Y=\PP^2\times \PP^2\times \PP^2$ in the case \eqref{list-rho=3-degree=30}, and
$Y=\PP^1\times \PP^1\times \PP^1\times \PP^1$ in the case \eqref{list-P1P1P1P1}.
The embedding $X \subset Y$ is canonical. 
This means, in particular, that automorphisms of $X$ are induced by 
those of $Y$.
\end{case}

Below we give some comments on the varieties in the table.

\begin{case}{\bf Case \eqref{list-rho=2-degree=12}.}
The variety $X$ can be one of the following two types:  
\begin{enumerate}
 \item[a)]
$X\subset \PP^2\times \PP^2$ is a divisor of 
 bidegree $(2,2)$, or
\par\medskip
 \item[b)]
 $X$ is a double cover of $V_6\subset \PP^2\times \PP^2$ (see \eqref{list-W6}),
 whose branch locus is a member of $|-K_{V_6}|$.
\end{enumerate}
\end{case}

\begin{case}{\bf Case \eqref{list-section-P3P3}.}\label{remark-list-section-P3P3}
The variety $X$ is the blow-up of $\PP^3$ along a 
 curve of degree $6$ and genus $3$ which is an intersection of cubics.
Indeed, two projections $f_i: X\to \PP^3$ are blowups of some smooth curve $C\subset \PP^3$.
Easy computations show that $C$ must satisfy the above conditions.

Conversely, let $X$ be  the blow-up of $\PP^3$ along a 
 curve of degree $6$ and genus $3$ which is an intersection of cubics.
By \cite[\S 5]{Mori1983} on $X$ there are two blowup structures
$f_i: X\to \PP^3$ of this type.
These two contractions induce a map 
$f=f_1\times f_2: X\to \PP^3\times \PP^3$ which must be finite and birational 
onto its image. 
Consider the composition $f': X \stackrel{f}{\to}\PP^3\times \PP^3\stackrel{s}{\to} \PP^{15}$,
where $s$ is the Segre embedding. Let $M_i:=f_i^* \OOO_{\PP^3}(1)$.
Then $-K_{X}=M_1+M_2$. By Riemann-Roch $\dim |-K_X|=12$. 
Therefore, $f'(X)$ is contained into a subspace of codimension $3$, i.e.,
$f'(X)\subset s(\PP^3\times \PP^3)\cap \PP^{12}$. 
Since $\deg f'(X) =20=\deg s(\PP^3\times \PP^3)\cap \PP^{12}$,
we have $f(X)\simeq f'(X)= s(\PP^3\times \PP^3)\cap \PP^{12}$.
By subadjunction, $K_X=f^*K_{f(X)}-B$, where $B$ is an effective divisor 
defined by the conductor ideal. Since, $\omega_{f(X)}=\OOO_{f(X)}(-1,-1)$,
we have $B=0$. Therefore, $f(X)$ is normal and $f(X)\simeq X$.

\begin{scase}\label{example-cubo-cubic-Cremona}
{\bf Example (cf. \cite{Katz1987}).}
Let $(a_{i,j})$, $(b_{i,j})$, $(c_{i,j})$ be 
symmetric $4\times 4$-matrices and let $X\subset \PP^3_{x_1\dots,x_4}\times \PP^3_{y_1\dots,y_4}$
is given by the equations
\[
\sum a_{i,j} x_i y_j=\sum b_{i,j} x_i y_j=\sum c_{i,j} x_i y_j=0.
\]
If $(a_{i,j})$, $(b_{i,j})$, $(c_{i,j})$ are taken sufficiently general, then
$X$ is a smooth Fano threefold of type \eqref{list-section-P3P3}.
The group $\muu_2$ acts on $X$ by $x_i \mapsto y_i$.
The induced birational involution $\PP^3\dashrightarrow \PP^3$ 
is cubo-cubic, it is given by the linear system
of cubics passing through the center of blowup $C\subset \PP^3$.
The exceptional divisor is sweep out by trisecants of $C$
and it is a surface $F\subset \PP^3$ of degree $8$ 
with multiplicity $3$ along $C$ \cite{Katz1987}.
\end{scase}\end{case}

\begin{case}{\bf Case \eqref{list-2-blowup-Q}.}
It is easy to show that $Y'$ is a Fano fivefold.
Linear systems $|\OOO_{\PP^5}(1)|$ and $|\OOO_{\PP^5}(2)-E|$ define 
two contractions $h_i: Y'\to \PP^5$ which are blowups of Veronese surfaces
$V=V_4\subset \PP^5$  (cf. \cite{Crauder1989}).
In this case,
$X$ can be realized as the blow-up of a smooth quadric $Q\subset \PP^4$ along a 
 twisted quartic curve.
Indeed, the restriction $f_1=h_1|_X$ is a birational map whose image is 
a quadric  $Q = h_1(D_1)\cap h_2(D_2) \subset \PP^5$ and moreover $f_i$ is the blowup of 
$Q$ along $V\cap h_1(D_1)$, a twisted quartic curve.
\begin{scase}
{\bf Example.} 
Let $C \subset \PP^4$ be a rational normal curve of degree $4$.
The action of the group $\Aut(C)=\PGL_2(\Bbbk)$ naturally extends to $\PP^4$
so that $\PP^4=\PP(S^4V)$, where $V$ is the standard representation of $\GL_2(\Bbbk)$.
The representation of $\GL_2(\Bbbk)$ on $S^4V$ is 
irreducible and can be defined by matrices over $\RR$.
Therefore, there exists an invariant non-singular quadric $Q\subset \PP^4$
containing $C$.
Let $f_1:X\to Q$ be the blowup of $C$. Then 
$X$ is a Fano threefold of type \eqref{list-2-blowup-Q}. 
The $\PGL_2(\Bbbk)$-action 
lifts to $X$ and the second contraction $f_2: X\to Q$ 
is also $\PGL_2(\Bbbk)$-equivariant. Let $C'\subset Q$ 
be the center of the blowup $f_2$.
Clearly, the pairs $(Q,C)$ and $(Q,C')$ are $\SL_2(\Bbbk)$-isomorphic, i.e. there exists 
an automorphism 
$\gamma: Q\to Q$ such that $\gamma(C)=C'$ (because $C$ is the only 
one-dimensional orbit on $Q$).
Then $f_2\comp \gamma \comp f_1^{-1}$ is an involution on $X$.
\end{scase}\end{case}

\begin{case}{\bf Case \eqref{list-rho=3-degree=30}.}
Here  $X$ is the blow-up of $V_6\subset \PP^2\times \PP^2$ (see \eqref{list-W6}),
 along a 
 curve $C$ of bidegree $(2,2)$ 
 such that the composition 
 $C\hookrightarrow V_6 \hookrightarrow \PP^2\times \PP^2 \stackrel{p_i}{\rightarrow}\PP^2$
 is an embedding for each projection $p_i$, $i=1,2$.
Indeed, there are three projections $\pi:= X\to \PP^2\times \PP^2$.
The image of each $\pi_i$ is contained into a divisor $V\subset \PP^2\times \PP^2$
of type $(1,1)$ and $\pi_i$ is birational onto its image.
Hence $\pi_i$ passes through a birational extremal contraction $X\to V'\to V$.
By \cite[Table 3]{Mori1981-82} we have $V'\simeq V$, $V$ is smooth and 
$\pi_i$ is a blowup of a curve as above. 

\begin{scase}
{\bf Example (cf. \cite{Nakano-1989}).} 
Let $\Gamma \subset \PP^2$ is a non-degenerate conic  and let $\Gamma^*\subset \PP^{2*}$
be its dual, the conic formed by lines that are tangent to $\Gamma$.  Consider the incidence curve
\[
C=\{ (P,L) \in \Gamma\times \Gamma^*\subset \PP^2\times \PP^{2*} \mid \text{$L$ is tangent to $\Gamma$ at $P$}\}.
\]
Then $C$ is contained into the flag variety $\Fl(\PP^2)=V_6$ and satisfies conditions of \eqref{list-rho=3-degree=30}.
The action $\Aut(\Gamma)=\PGL_2(\Bbbk)$ extends to $V_6=\Fl(\PP^2)$.
Orbits of $\PGL_2(\Bbbk)$  on $V_6$ are described as follows:
\begin{itemize}
\item
$C\simeq \PP^1\simeq \SL_2(\Bbbk)/B$, where $B$ is a Borel subgroup,
\item 
$D'$ and $D''$, where $\bar D':= D'\cup C$ and $\bar D'':= D''\cup C$ 
are complete surfaces isomorphic to $\PP^1\times \PP^1$ with the diagonal action of $\PGL_2(\Bbbk)$,
\item
an open orbit $U\simeq \SL_2(\Bbbk)/Q_8$, where $Q_8$ is the binary 
quaternion group of order $8$.
\end{itemize}
Let $f: X\to V_6$ be the blowup of $C$. Then $X$ is a Fano threefold of type \eqref{list-rho=3-degree=30}
admitting a $\PGL_2(\Bbbk)$-action.
There are two more
$\PGL_2(\Bbbk)$-equivariant contractions $f', f'': X\to V_6$ contracting proper transforms 
of  $\bar D'$ (resp. $\bar D''$) 
to curves $C'\subset V_6$ (resp.  $C''\subset V_6$) of bidegree $(2,2)$. 
The pairs $(V_6, C)$,  $(V_6, C')$, and $(V_6, C'')$, are isomorphic.
These isomorphisms induce an action of the symmetric group  $\Sym _3$ on $X$.
\end{scase} \end{case} 

\begin{case}{\bf Case \eqref{list-P1P1P1P1}.}
The variety $X$ is isomorphic to the blowup of $\PP^1\times \PP^1\times \PP^1$ along 
an elliptic curve which is an intersection of two 
members of $|-\frac12K_{\PP^1\times \PP^1\times \PP^1}|$.
Similar to \ref{example-cubo-cubic-Cremona}
one can easily construct examples of symmetric 
varieties of this type (admitting a $G$-structure). 
\end{case}

\section{Action on the Picard lattice.}\label{section-Action-on-the-Picard-lattice}

\begin{mtheorem}{\bf Lemma.}
\label{lemma-representation1}
Let $V=\QQ^N$ and let $\Phi: G\hookrightarrow \GL(V)$ be a faithful 
representation of a finite group $G$.
Identify $G$ with its image $\Phi(G)\subset \GL(V)$ and assume that $V^G=0$.
Then there is a subgroup $G_0\subset G$ such that 
$V^{G_0}=0$ and
\[
|G_0|=
\begin{cases}
\text{$2$}&\text{if $N=1$,}
\\
\text{$2$ or $3$}&\text{if $N=2$,}
\\
\text{$2$, $4$, or $6$}&\text{if $N=3$,}
\\
\text{$5$, $9$, or divides $3\cdot 2^7$}&\text{if $N=4$.}
\end{cases}
\]
\end{mtheorem}
\begin{proof}
Cases $N\le 2$ are trivial. Consider the case $N=3$.
According to Minkowski's bound (see e.g. \cite{Serre2007})
the order of any subgroup $G\subset \GL_3(\QQ)$ divides 
$48$.
Assume that there is an element $\tau\in G$ such that $V^{\langle\tau\rangle}=0$.
Let $\mu(t)$ be 
the minimal polynomial of $\tau$.
Then $\deg \mu\le 3$, $\mu$ is a product of cyclotomic polynomials $\phi_k(t)$,
and $\mu(1)\neq 0$.
Moreover, in some basis, $\tau$ is given by an orthogonal matrix,
so $\mu(-1)=0$ and $\phi_2 \mid \mu$. 
Hence, there are the following possibilities:
$\mu=\phi_2$, $\phi_2\phi_4$, $\phi_2\phi_3$, or $\phi_2\phi_6$.
Thus we can take $G_0=\langle\tau\rangle$.

Now we assume that for any element $\tau\in G$ we have $V^{\langle\tau\rangle}\neq 0$. 
In particular, $G$ is not a cyclic group. 
We also may assume that $|G|\ge 8$.
If $G$ contains a subgroup $G_1\subset G$ isomorphic to
$\muu_2\times \muu_2$, then either 
$G_1\subset \SL(V)$ or $G_1\ni -\id$.
In the latter case we can take $G_0=\langle -\id\rangle$ and
in the former case we can take $G_0=G_1$.
Note that $\GL(V)$ contains no elements of order $8$
and the quaternion group has no real
faithful 
representations of dimension $\le 3$.
Thus we may assume that the order of $G$ is not divisible by $8$.
We get only one possibility:
Sylow $2$-subgroups of $G$ are cyclic of order $4$
and $G$ is of order $12$.
Let $G_1\subset G$ be a Sylow $2$-subgroup
and let $\tau$ be its generator.
We may assume that $G_1\subset \SL(V)$
(otherwise we can take $G_0=G_1$).
Then $G\subset \SL(V)$. If $G$ has an element $\xi$ of order $6$, then $\xi$
generates the desired subgroup $G_0$.
Otherwise the representation $\Phi: G\hookrightarrow \SL(V)$
is irreducible. Then we must have $G\simeq \mathfrak{A}_4$
and $G_1$ is not cyclic, a contradiction.

Consider the case $N=4$. Again according to Minkowski's bound (see e.g. \cite{Serre2007})
the order of any subgroup $G\subset \GL_4(\QQ)$ divides $2^7 \cdot 3^2 \cdot 5$.
Again if there is an element $\tau\in G$ such that $V^{\langle\tau\rangle}=0$, then
its minimal polynomial $\mu(t)$ satisfies the following conditions: 
$\deg \mu\le 4$, $\mu$ is a product of cyclotomic polynomials $\phi_k(t)$,
and $\mu(1)\neq 0$. In this case, we have either $\mu=\phi_5$, $\phi_8$, $\phi_{10}$, $\phi_{12}$,
or $\mu$ divides $\phi_2\phi_3\phi_4\phi_6$.
Thus we can take $G_0=\langle\tau\rangle$.

Now we assume that for any element $\tau\in G$ we have $V^{\langle\tau\rangle}\neq 0$. 
Then $G$ contains no elements of order $5$, $9$, and $8$. 
In particular, the order of $G$ divides $2^7 \cdot 3^2$.
Let $G_1\subset G$ be a Sylow $3$-subgroup. By the above, we may assume that 
$G_1\simeq \muu_3\times \muu_3$. In this case $V^{G_1} = \{0\}$
and we can take $G_0=G_1$.
\end{proof}

\begin{stheorem}{\bf Corollary.}
\label{corollary-representation}
In notation of Theorem \xref{main}
there is a subgroup $G_0\subset G$ of order $N$ such that 
$\rho(X)^{G_0}=1$, where 
\[
N=
\begin{cases}
2&\text{if $\rho(X)=2$,}
\\
\text{$2$ or $3$}&\text{if $\rho(X)=3$,}
\\
\text{$2$, $4$ or $6$}&\text{if $\rho(X)=4$,}
\\
\text{$5$, $9$, or divides $3\cdot 2^7$}&\text{if $\rho(X)=5$.}
\end{cases} 
\]
\end{stheorem}

\begin{proof}
Let $V$ be the orthogonal complement to $c_1$ in $\Pic_{\QQ}(X):=\Pic(X)\otimes \QQ$. 
The faithful representation $G\hookrightarrow \GL(V)$
satisfies the condition $V^G=0$. Then we can apply Lemma \ref{lemma-representation1}.
\end{proof}

\section{Fano conic bundles and del Pezzo fibrations.}
\label{section-Fano-conic-bundles-and-del-Pezzo-fibrations}

\begin{case}{\bf Definition.}
Let $X$ be a smooth threefold. A morphism $f:X\to Z$ onto a smooth surface 
is a \textit{conic bundle} if every (scheme) fiber is isomorphic to a conic 
in $\PP^2$. A morphism $f:X\to Z$ onto a smooth curve 
is a \textit{del Pezzo bundle} if $f$ has connected fibers and $-K_X$ is $f$-ample.
(In this case a general fiber of $f$ is a del Pezzo surface.)
\end{case}

\begin{scase}{\bf Remark.}
Let $X$ be a smooth Fano threefold and let $f: X\to Z$ 
be a surjective morphism with connected fibers.
\begin{enumerate}
\item 
If $Z$ is a smooth curve, then  $f$ is a del Pezzo bundle
and $Z\simeq \PP^1$.
\item 
If $Z$ is a smooth surface and $f$ is equidimensional, then $f$ is a 
conic bundle and $Z$ is rational.
\end{enumerate}
\end{scase}

\begin{stheorem}{\bf Proposition (\cite[\S 4]{Mori-Mukai-1986}).}
\label{Proposition-Fano-conic-bundles}
Let $X$ be a Fano threefold. 
Assume that $X$ has a conic bundle structure $f: X\to Z$. Then
\begin{enumerate}
 \item 
$Z$ is a del Pezzo surface.
 \item 
The discriminant curve $\Delta\subset Z$ is a curve with at worst ordinary double points
\textup(or empty\textup).
 \item 
$\rho(X/Z)=1$ if and only if for any  irreducible curve $C\subset Z$ its 
preimage $f^{-1}(C)$ is also irreducible.
 \item 
If $C\subset Z$ is a an irreducible curve such that $f^{-1}(C)$ is reducible, 
then $C$ is a smooth connected component of $\Delta$.
 \item 
$h^{1,2}(X)=\rho(X)-\rho(Z) +p_a(\Delta)-2$.
\end{enumerate}
\end{stheorem}

\begin{stheorem}{\bf Corollary.}
If in the assumptions of Proposition \xref{Proposition-Fano-conic-bundles}
$\rho(Z)\ge 2$, then $X$ has a del Pezzo bundle structure.
\end{stheorem}

\begin{case}
{\bf Assumption.}
>From now on we assume that $X$ is a smooth Fano threefold 
satisfying assumptions of Theorem \xref {main}.
By \cite{Prokhorov2010a} we may assume that 
the Fano index of $X$ is equal to $1$ (otherwise $X$ is of type 
\eqref{list-W6} or \eqref{list-P1P1P1}).
Thus from now on we assume that $\Pic(X)^G$ is generated by $-K_X$. 
We also assume that $G$ is the smallest group satisfying  condition \eqref{condition}
(cf. Corollary \xref{corollary-representation}).
\end{case}

\begin{mtheorem}{\bf Proposition.}\label{conic-bundle}
Assume that $X$ has a conic bundle structure 
$f:X\to Z$ over $Z=\PP^2$.  Then $X$ is of type \eqref{list-rho=2-degree=12} or \eqref{list-rho=3-degree=30}.
\end{mtheorem}

\begin{proof}
According to \cite[Prop. 6.3]{Mori1983}
$\rho(X)\le 3$.
By Corollary \ref{corollary-representation} we may assume that $|G|\le 3$.
Let $L\subset \PP^2$ be a line, 
let $F:=f^{-1}(L)$, and let 
$\Delta\subset \PP^2$ be the discriminant curve.
Let $F_1,\dots, F_n$ be the $G$-orbit of 
the class of $F$ in $\Pic(X)$, where $n=2$ or $3$. Write $\sum F_i=ac_1$ for some $a\in \ZZ$. 
Take the line $L$ to be sufficiently general, so the surface $F$ is smooth. 
Then we have 
\begin{multline}\label{eq-disc-curve}
n(12-\deg \Delta)=nK_F^2+4n =n(K_X+F)^2\cdot F-
\\
-2nc_1\cdot F^2=nc_1^2\cdot F=c_1^2\cdot \sum F_i=ac_1^3>0.
\end{multline}
In particular, $\deg \Delta<12$.

First assume that $n=2$.
As above $F_1+F_2=ac_1$ and so 
\[
a^2c_1^3= \frac 1a (F_1+F_2)^3=\frac 6a F_1^2\cdot F_2
=\frac 6a F_1^2\cdot (F_1+F_2)=6c_1\cdot F^2=12.
\]
Since $c_1^3$ is even, 
there is only one possibility:
$a=1$, $c_1^3=12$, and $\deg \Delta=6$. 
By Proposition \ref{Proposition-Fano-conic-bundles}
$h^{1,2}(X)=\rho(X)+7$.
Then from tables in \cite{Mori1981-82} we get the
 case \eqref{list-rho=2-degree=12}.
 
Assume that $n=3$. Then  $\rho(X)=3$ and $\Delta\neq \emptyset$.  From \eqref{eq-disc-curve} 
we see $3(12-\deg \Delta)=ac_1^3$.
So, $ac_1^3\le 33$.

If $a\ge2$, then from \cite[Table 3]{Mori1981-82}
we get $a=2$. Then $c_1^3=18-\frac32 \deg \Delta$.
Since $c_1^3$ is even, $c_1^3=12$ and $\deg \Delta=4$.
But in this case, $h^{1,2}(X)=p_a(\Delta)=3$.
This contradicts \cite[Table 3]{Mori1981-82}.
Therefore, $a=1$ and $c_1^3=3(12-\deg \Delta)$.
In particular, $\deg \Delta$ is even and $c_1^3\le 30$.
Moreover, if $c_1^3=30$, then by 
\cite[Table 3]{Mori1981-82}
we get the case \eqref{list-rho=3-degree=30}.
Thus we may assume that $c_1^3<30$.
There are the following possibilities:
\begin{itemize}
\item 
$\deg \Delta=4$, $c_1^3=24$, $h^{1,2}(X)=p_a(\Delta)=3$,
\item 
$\deg \Delta=6$, $c_1^3=18$, $h^{1,2}(X)=p_a(\Delta)=10$,
\item 
$\deg \Delta=8$, $c_1^3=12$, $h^{1,2}(X)=p_a(\Delta)=21$.
\end{itemize}
All these cases are impossible by the classification \cite[Table 3]{Mori1981-82}.
\end{proof}

\begin{mtheorem}{\bf Lemma.}
\label{del-Pezzo-bundle}
Assume that $X$ has a del Pezzo bundle structure $f:X\to \PP^1$ and let 
$F$ be a general fiber. Let $F_1,\dots, F_n$ be the $G$-orbit of 
the class of $F$ in $\Pic(X)$. Write $\sum F_i=ac_1$ for some $a\in \ZZ$. 
Then 

\begin{enumerate}
\item
$9n\ge nK_F^2=ac_1^3$, \ $a>0$,
\item
$n>2$, $|G|>2$, and $\rho(X)\ge 3$,
\item
$ac_1^3$ is divisible by $3$,  
\item
$aK_F^2$ is even,
\item
if $n=3$, then $\rho(X)= 3$ and either 
\begin{enumerate}
\item
 $X$ is of type \eqref{list-double-cover-P1P1P1},
or 
\item
$c_1^3=24$ and $K_F^2=8$.
\end{enumerate}
\end{enumerate}
\end{mtheorem}

\begin{proof}
We have 
\begin{equation*}
9n\ge nK_F^2=nc_1^2\cdot F=c_1^2\cdot \sum F_i=ac_1^3.
\end{equation*}
This proves (i). Further,
\begin{equation}
6\sum_{i<j<k} F_i\cdot F_j\cdot F_k=\bigl(\sum F_i\bigr)^3 =a^3c_1^3>0.
\end{equation}
In particular, $ac_1^3$ is divisible by $3$ and $n>2$.
Further,
\[
2F\cdot \sum_{1<i<j} F_i\cdot F_j=F\cdot \bigl(\sum F_i\bigr)^2=a^2F\cdot c_1^2=a^2K_F^2. 
\]
Hence, $a^2K_F^2$ is even.

Finally, let $n=3$. Then $\rho(X)\ge 3$.
If $c_1^3\le 12$, then by \cite[Tables 3-5]{Mori1981-82} we are in the 
case \eqref{list-double-cover-P1P1P1}.
Thus we may assume that $c_1^3>12$.
Since $c_1^3$ is even, 
$ac_1^3\le 24$. Therefore, $a=1$.
Let $\Gamma:=F_2|_F$. Further,
\[
-K_F\cdot \Gamma=c_1\cdot F_2\cdot F_1= 
F_1\cdot F_2\cdot F_3=\frac16(F_1+F_2+F_3)^3=\frac{c_1^3}{6}.
\]
Since $\Gamma^2=F_2^2\cdot F=0$,
by Riemann-Roch $K_F\cdot \Gamma$ is even. Hence 
$c_1^3$ is divisible by $12$.
We get $c_1^3=24$ and $K_F^2=8$. 
If $\rho(X)\ge 4$, then
by \cite[Tables 4-5]{Mori1981-82} we have $\rho(X)= 4$ and $X$ is of type \eqref{list-P1P1P1P1}.
But in this case, $X$ has a del Pezzo bundle structure of degree $6$, a contradiction.
\end{proof}

\begin{mtheorem}{\bf Lemma.}\label{lemma-product}
$X\not \simeq Z\times \PP^1$, where $Z$ is a smooth surface. 
\end{mtheorem}
\begin{proof}
Clearly, $Z$ is a del Pezzo surface of degree $10-\rho(Z)$.
Let $F$ be a fiber of the projection $X=Z\times \PP^1\to \PP^1$.
Take an element $\tau\in G$ so that ${}^\tau F$ is not proportional to $F$.
Then ${}^\tau F\sim \alpha F+ f^* L$ for some $0\neq L\in \Pic(Z)$
and $\alpha \in \ZZ$.
Since $F^2\equiv 0$, we have 
\[
0= {}^\tau F^2\cdot F= f^*L^2\cdot F.
\]
Hence, $L^2= 0$ and $2 \alpha F\cdot f^* L\equiv {}^\tau F^2 \equiv 0$. 
So, $\alpha=0$ and ${}^\tau F= f^* L$. Further, by Riemann-Roch $K_Z\cdot L$
is even and
\[
K_Z^2=c_1^2\cdot F=c_1^2\cdot {}^\tau F=c_1^2\cdot f^* L=
(2F-f^*K_Z)^2\cdot f^* L=-4 K_Z\cdot L.
\]
Therefore, $K_Z^2=8$ and $\rho(X)=3$. This contradicts
Lemma  \ref{del-Pezzo-bundle}. 
\end{proof}

\section{Proof of Theorem \ref{main}}\label{section-Proof-of-Theorem}
Recall our assumption that $\Pic(X)^G$ is generated by $c_1$ and that $G$ is the smallest group satisfying  condition \eqref{condition}.
\begin{case}
Consider the case $\rho(X)=2$.
Denote the generator of $G\simeq\muu_2$ by
$\tau$.
Let $R_1$ and $R_2$ be extremal rays of the Mori cone
$\NE(X)\subset \N_1(X)\simeq \RR^2$.
Let $f_i:X\to X_i$ be the contraction of $R_i$, let 
$A_i$ be the ample generator of $\Pic(X_i)\simeq\ZZ$,
and let $M_i:=f_i^*A_i$.
By Proposition \ref{conic-bundle} and Lemma \ref{del-Pezzo-bundle}
we may assume that both contractions are birational.
Let $D_i$ be the exceptional divisor of $f_i$.

\begin{stheorem}{\bf Theorem (\cite[Th. 5.1]{Mori1983}).}
The group $\Pic(X)$ is generated by $M_1$ and $M_2$.
\end{stheorem}
By this theorem 
\begin{equation*}
\label{rho2-eq-K}
-K_X\equiv \alpha_1M_1+\alpha_2M_2,\quad \alpha_i\in \ZZ
\end{equation*}
Since $M_1+{}^\tau M_1$ and $M_2+{}^\tau M_2$ are invariant divisors,
we have 
\[
M_1+{}^\tau M_1=-a_1K_X,\qquad 
M_2+{}^\tau M_2=-a_2K_X,
\]
for some $a_1, \, a_2\in \ZZ$. Then
\[
 -2K_X=\alpha_1M_1+\alpha_2M_2+{}^\tau(\alpha_1M_1+\alpha_2M_2)=-(\alpha_1a_1+\alpha_2a_2)K_X.
\]
So, $\alpha_1a_1+\alpha_2a_2=2$. On the other hand,
\[
 0< 2 K_X^2\cdot M_i= K_X^2\cdot (M_i+{}^\tau M_i)=a_ic_1^3.
\]
Therefore, $a_i>0$.
This gives us $a_i=\alpha_i=1$, i.e., 
\begin{equation}
\label{rho2-eq-K-1}
{}^ \tau M_1 = M_2,\qquad -K_X=M_1+M_2.
\end{equation}
Further, $M_1^2\cdot D_1=M_2^2\cdot D_2=0$ and   $M_2^2\cdot {}^ \tau D_1=  {}^ \tau M_1^2\cdot {}^ \tau D_1=0$.
Hence, $D_2$ and ${}^ \tau D_1$ are proportional, so ${}^ \tau D_1=bD_2$ for some $b\in \QQ$.
Since  $D_2\cdot (-K_X)^2>0$ and   ${}^ \tau D_1\cdot (-K_X)^2>0$, $b>0$.
Thus 
\[
D_1+bD_2=D_1+{}^ \tau D_1=-cK_X \quad \text{for some $c\in \ZZ_{>0}$}.
\]

If $f_1(D_1)$ is a point, then  $M_1\cdot D_1\equiv M_2\cdot D_2\equiv 0$, i.e.
$f_2(D_2)$ is also a point. In this case, $D_1\cap D_2=\emptyset$.
Hence, 
\[
0<c^2(-K_X)^3=(-K_X)\cdot (D_1+bD_2)^2=(-K_X)\cdot D_1^2+b^2(-K_X)\cdot D_2^2<0,
\]
a contradiction.

Therefore, $C_i:=f_i(D_i)$ are curves.
In this case both varieties $X_i$ are smooth Fano threefolds. Write 
$-K_{X_i}=r_i A_i$, $r_i=1$, $2$, $3$, or $4$.
Since $D_2$ and   ${}^ \tau D_1$ are primitive elements of $\Pic(X)$,
we have $D_2={}^ \tau D_1$.
Then 
\[
-K_X=-f^*K_{X_1}-D_1=r_1M_1-D_1 = r_1{}^ \tau M_1-{}^ \tau D_1= r_1M_2-D_2.
\]
This shows that ${}^\tau D_1=D_2$ and $r_1=r_2$. Put $r:=r_1=r_2$.
Further, 
\[
M_2^2\cdot M_1=M_1^2\cdot M_2=M_1^2\cdot ( (r-1)M_1-D_1)=(r-1)M_1^3,
\]
\[
A_1\cdot C_1=-M_1\cdot D_1^2= -M_1\cdot ((r-1)M_1-M_2 )^2=
(r-1)(r-2)M_1^3.
\]
Therefore, $r\ge 3$. 
We get two possibilities: \eqref{list-section-P3P3} and \eqref{list-2-blowup-Q}.
\end{case}

>From now on we assume that $\rho(X)\ge 3$.

\begin{mtheorem}{\bf Proposition.}\label{proposition-del_pezzo-bundle-rho3}
If in the above assumptions $\rho(X)\ge 3$ and $X$ is not of type \eqref{list-rho=3-degree=30}, then $X$
has a structure of del Pezzo bundle of degree $\le 8$.
\end{mtheorem}
\begin{proof}
If $X$ has a conic bundle structure $f: X\to Z$, then by 
Proposition \ref{conic-bundle} $Z\not\simeq\PP^2$.  Since $Z$ is a 
smooth rational surface, there exists a surjective morphism
$Z\to \PP^1$ with connected fibers. The composition map gives 
a structure of del Pezzo bundle whose fibers $F$ are del Pezzos 
with $\rho(F)>1$. 

Assume that $X$ has no conic bundle structures.
Then by \cite[(9.1), (9.2), (9.6))]{Mori1983} the variety
$X$ isomorphic to the blow-up of $\PP^3$ along a disjoint union of 
a line and a conic.
The map $X\to \PP^3$ can be decomposed 
$X\to X_1\to\PP^3$, where $X_1\to\PP^3$ is the blow-up of a line.
In this case $X_1$ has a structure of a $\PP^2$-fibration over $\PP^1$.
The composition $X\to X_1\to\PP^1$ is the desired map.
\end{proof}

\begin{case}
Assume that $\rho(X)=3$.
By Proposition \ref{proposition-del_pezzo-bundle-rho3} 
$X$ has a del Pezzo bundle structure $f:X\to Z$
and by Lemma \ref{del-Pezzo-bundle}  $n:=|G|>2$.
Let $F$ be a general fiber.
By Lemma \ref{del-Pezzo-bundle} we have $K_F^2=8$ and $c_1^3=24$.
According to \cite[Table 3]{Mori1981-82} we have the following possibilities:

\begin{scase} {\bf Case $7^o$.}
$X$ is a blow-up of $V_6\subset \PP^2\times \PP^2$
along an elliptic curve $C$ which is an intersection of two
members of $|-\frac 12 K_{V_6}|$.
Then projection from $C$ gives us 
a del Pezzo bundle structure $f: X\to \PP^1$
of degree $6$.
This contradicts 
$K_F^2=8$.
\end{scase}

\begin{scase}
{\bf Case $8^o$.}
$X$ is a member of the linear system $|p_1^*\sigma^*\OOO(1)\otimes p_2^*\OOO(2)|$
on $\FF_1\times \PP^2$, where $p_i$ is the projection and $\sigma:\FF_1\to \PP^2$
is the blowing up. The composition $X\stackrel {p_1}{\longrightarrow} \FF_1\stackrel {q}{\longrightarrow}\PP^1$,
where $q$ is the $\PP^1$-ruling, is a del Pezzo bundle whose general fiber $F$ is a 
divisor in $\PP^1\times \PP^2$ of bidegree $(1,2)$. 
Hence, $F$ is a del Pezzo surface of degree $5$.
This contradicts 
$K_F^2=8$.
\end{scase}
\end{case}

\begin{case}
Finally, let $X$ be a Fano threefold with $\rho(X)\ge 4$.
By \cite[Tables 4-5]{Mori1981-82} $c_1^3\ge 24$.
Moreover, if $c_1^3= 24$, then $X$ is a divisor of multidegree
$(1,1,1,1)$ in $\PP^1\times \PP^1\times \PP^1\times \PP^1$, i.e. $X$ is of type 
\eqref{list-P1P1P1P1}. 
Thus from now on we assume that $c_1^3\ge 26$
and $4\le \rho(X)\le 5$
(see \cite[Tables 4-5]{Mori1981-82}).
\end{case}

\begin{case}
Assume that $\rho(X)=4$. Then by Corollary \ref{corollary-representation} $n=2$, $3$, $4$, or $6$.
According to Lemma \ref{del-Pezzo-bundle} (i)
\begin{equation}
\label{eq-rho=4}
48\ge nK_F^2=ac_1^3
\end{equation}
and  $ac_1^3$ is divisible by $3$.

Since, $c_1^3\ge 26$, 
we have $a=1$, $nK_F^2\ge 26$, and $n=4$ or $6$. 
Further, by Lemma  \ref{del-Pezzo-bundle} (iv)
$K_F^2$ is even. Hence, $c_1^3$ is divisible by $12$. 
By \cite[Table 4]{Mori1981-82} we get the following possibility
$7^o$:
$X$ is the blow-up 
of $V_6\subset \PP^2\times \PP^2$ of type \eqref{list-W6} along 
disjoint union of curves of bidegree $(1,0)$ and $(0,1)$.
So, there is an embedding $X\subset \FF_1\times \FF_1$.
The composition
\[
X\hookrightarrow \FF_1\times \FF_1 \stackrel{p_1}{\longrightarrow}
 \FF_1 \to \PP^1
\]
is a del Pezzo bundle whose general fiber 
is a del Pezzo surface of degree $7$.
This contradicts \eqref{eq-rho=4}.
\end{case}

\begin{case}
Assume that $\rho(X)=5$. 
Then by \cite[Table 5]{Mori1981-82} and Lemma \ref{lemma-product} we have two possibilities:

\begin{scase}
\textbf{Case $1^o$.}
Then 
$c_1^3=28$ and $X$ can be obtained as a 
sequence of blow-ups
\[
X\stackrel{h}{\longrightarrow} Y \stackrel{g}{\longrightarrow} Q,
\]
where $Q$ is a smooth quadric in $\PP^4$, $g$ is the blow-up of a conic 
$C\subset Q$, and $h$ is the blow-up of three exceptional 
lines of $g$.
As above, the projection $Q\dashrightarrow \PP^1$ from $C$ induces 
a quadric bundle $Y\to \PP^1$. The composition 
$s: X\stackrel{h}{\longrightarrow} Y \to \PP^1$ 
is a del Pezzo bundle whose general fiber $F$ is 
the blow-up of three points on the corresponding fiber of
$Y\to \PP^1$. Hence, $F$ is a del Pezzo surface of degree $5$.
By Lemma \ref{del-Pezzo-bundle} (i) we have $5n=nK_F^2= 28a$.
Since $|G_0|$ is not divisible by $7$, we get a contradiction.
\end{scase}

\begin{scase}
\textbf{Case $2^o$.}
$c_1^3=36$ and $X$ can be obtained as a 
sequence of blow-ups
\[
X\stackrel{h}{\longrightarrow} Y \stackrel{g}{\longrightarrow} \PP^3,
\]
where $g$ is the blow-up of two disjoint lines $L_1,\, L_2\subset \PP^3$ 
and $h$ is the blow-up of two exceptional 
lines $\ell,\, \ell' \subset g^{-1}(L_1)$.
The projection from $L_2$ induces a del Pezzo bundle 
$s: X\to \PP^1$ of degree $8$ and projection from $L_1$ induces a del Pezzo bundle 
$r: X\to \PP^1$ of degree $6$.
Let $F$ be a general fiber of $s$. Then by Lemma \ref{del-Pezzo-bundle} (i)
\[
8n=-a K_X^3= 36a, \quad 2n=9a.
\]
Hence $9$ divides $n$ and we may assume that $n=9$ and $a=2$.
By Lemma \ref{lemma-representation1} and our assumption $|G|=9$.
Let $F'$ be a general fiber of $r$. Again by Lemma \ref{del-Pezzo-bundle} (i)
\[
6n'=-a' K_X^3= 36a', \qquad n'=6a'.
\]
On the other hand, $n'$ divides $|G|=9$, a contradiction.
\end{scase}
\end{case}
This finishes the proof of Theorem \ref{main}.

\section{Singular case}\label{section-Singular-case}
In this section we consider the case of singular $G$-Fano threefolds with $\rho(X)>1$.
On this step we may assume that the ground field is $\CC$.
\begin{mtheorem}{\bf Theorem (\cite{Namikawa-1997}).}
\label{theorem-Namikawa}
Let $X$ be a Fano threefold with terminal Gorenstein
singularities. Then $X$ is smoothable, that is,
there exists a flat family $\mathfrak X \to \mathfrak{D}\ni 0$ over a small disc $(\mathfrak{D}\ni 0)\subset \CC$ 
such that $\mathfrak X_0\simeq X$ and a general member $\mathfrak X_s$, $s\in \mathfrak{D}$ is 
a smooth Fano threefold.
Moreover, there is a natural identification 
$\Pic(X)=\Pic(\mathfrak X_s)=\Pic(\mathfrak X)$ so that 
$K_{\mathfrak X_s}=K_X$ \textup(see {\cite[\S 1]{jahnke-Radloff-2006arx}}\textup). 
\end{mtheorem}

\begin{case}\label{definition-contractible-plane}
Now let $X$ be a Fano threefold with terminal Gorenstein
singularities and let $\mathfrak X \to \mathfrak{D}\ni 0$ be its smoothing as above.
We say that a reduced  irreducible surface $S\subset X$ is a \textit{plane} if  
$S\simeq \PP^2$ and $\OOO_S(-K_X)=\OOO_{\PP^2}(1)$.
We say that a plane $S\subset X$ is \textit{contractible} if there exists a birational 
morphism $f: X\to Y$ to a normal variety $Y$ such that $\rho(X/Y)=1$,
$S$ is contained into the exceptional locus $\Exc(f)$ and 
$S$ does not meet other components of $\Exc(f)$ (if there is any).
Note that, in this situation, $\Exc(f)$ is of pure dimension $2$
(cf. \cite[Prop. 1.4]{Kachi1998}).

Recall that the \textit{nef cone} $\Nef(X)\subset H^2(X,\RR)$ is the closed cone generated by
nef divisors.
\end{case}

\begin{mtheorem}{\bf Proposition.}\label{Propositiontheorem-nef-cone}
Assume that, in the above notation, $X$ does not contain a contractible plane. Then
under  our identification
$\Pic(X)=\Pic(\mathfrak X/ \mathfrak D)=\Pic(\mathfrak X_s)$, $s\in \mathfrak{D}$  we have 
$\Nef(X)=\Nef(\mathfrak X/ \mathfrak D)=\Nef(\mathfrak X_s)$.
\end{mtheorem}

\begin{proof}
Note that $\mathfrak X$ has only isolated hypersurface singularities.
Hence, by \cite[Exp. XI, Corollary 3.14]{Grothendieck1968}
 the variety $\mathfrak X$ is (locally) factorial.
Let $\mathcal M$ be a divisor on $\mathfrak X$, let $M:=\mathcal M|_X$, and
$\mathcal M_s:=\mathcal M|_{\mathfrak X_s}$.
If $M$ is nef, then, obviously,  so $\mathcal M_s$ is. Thus it is sufficient to 
show the inverse implication. So, we assume that $\mathcal M_s$ is nef for 
some $s\in \mathfrak{D}$ and $M$ is not nef. Then $\mathcal M$ is also not nef.
Let 
\[
\lambda_0:= \sup \{ \lambda\in \QQ \mid \text{$\lambda\mathcal M-  K_{\mathfrak X}$ is nef over $\mathfrak{D}$} \}. 
\]
By the rationality theorem we have $\lambda_0\in \QQ$ and there is an extremal ray 
$R$ on $\mathfrak X/\mathfrak{D}$ such that $(\lambda_0\mathcal M-  K_{\mathfrak X})\cdot R=0$ and $ K_{\mathfrak X}\cdot R<0$.
Clearly,  $\lambda_0\mathcal M-  K_{\mathfrak X}$ is ample on $\mathfrak X_s$.
Therefore, the locus $\Exc(R)$ of the ray $R$ does not meet $\mathfrak X_s$, so
$\Exc(R)\subset X$. In particular, this means that $R$ is a flipping 
extremal ray. Then by  \cite{Kachi1998} the locus $\Exc(R)$ is a 
disjoint union of irreducible surfaces $S_i$ isomorphic to $\PP^2$
and, moreover, $\OOO_{S_i}(-K_{\mathfrak X})=\OOO_{\PP^2}(1)$.
Hence, $X$ contains a contractible plane, a contradiction.
\end{proof}

\begin{stheorem}{\bf Corollary.}\label{Corollary-nef-cone}
Notation as in \xref{definition-contractible-plane}.
Assume that $X$ does not contain a contractible plane.
Let $\mathfrak f_s: \mathfrak X_s\to \mathfrak Z_s$ be an extremal contraction.  
Then there exists an extremal contraction $\mathfrak f: \mathfrak X\to \mathfrak Z$
over $\mathfrak D$ such that the restriction $\mathfrak f|_{\mathfrak X_s}$ coincides with
$\mathfrak f_s$.
\end{stheorem}

\begin{mtheorem}{\bf Proposition.}\label{theorem-nef-cone_G}
Let $X$ be a $G$-Fano threefold. Then $X$ does not contain a contractible plane.
\end{mtheorem}

\begin{proof}
Assume the converse. Let $S\subset X$ be a contractible plane and let
$f : X\to Y$ be its contraction as in \ref{definition-contractible-plane}.
Let $A$ be an ample Cartier divisor on $Y$ and let $M:=f^*A$. 
Thus $M$ is a nef and big divisor on $X$ and $M$
is trivial on $S$. Now let $\{M_i\}$
and $\{S_{i,j}\}$ be $G$-orbits of $M$ and $S$, respectively, where 
the subscript index $i$ is chosen so that $M_i$ is trivial on $S_{i,j}$.
Then $M_i$ defines a contraction $f_i: X\to Y_i$ as in \ref{definition-contractible-plane} and 
$\Exc(f_i)=\cup_j S_{i,j}$.
For each fixed $i$, the surfaces $S_{i,j}$ do not meet each other.
Assume that the intersection $S_{i,j}\cap S_{i',j'}$ contains a curve $C$
for some $i\neq i'$. Then $M_i\cdot C=M_{i'}\cdot C=0$. So,  $C$ is contracted by
$f_i$ and $f_{i'}$. This contradicts $\rho(X/Y_i)=\rho(X/Y_{i'})=1$.
Therefore, any two different planes from $\{S_{i,j}\}$
intersect each other by a set of dimension $\le 0$.
On the other hand, $\sum_{i,j}S_{i,j}$ is a Cartier divisor proportional to 
$-K_X$. In particular, $\sum_{i,j}S_{i,j}$ is connected and  has 
only a local complete intersection singularities, a contradiction.
\end{proof}

Let $X$ be a $G$-Fano threefold with $\rho(X)>1$ (and  terminal Gorenstein singularities).
According to Theorem \ref{main} its smoothing 
is of one of the types \eqref{list-rho=2-degree=12}
--
\eqref{list-P1P1P1P1}. In this situation, we say that $X$ 
is of the corresponding type \eqref{list-rho=2-degree=12}
--
\eqref{list-P1P1P1P1}.

\begin{mtheorem}{\bf Theorem.}\label{theorem-not-QQ}
Let $X$ be a $G$-Fano threefold with $\rho(X)>1$.
\begin{enumerate}
 \item 
If $X$ is of type \eqref{list-W6} or \eqref{list-P1P1P1}, 
then $X$ is smooth.
 \item 
If $X$ is singular,
then $X$ has the same description as in 
the table of Theorem \xref{main}. 
 \item 
If $X$ is of type \eqref{list-2-blowup-Q} and $X$ is singular, then 
$X$ is the blowup of a quadratic cone in $\PP^4$
with center a union of two conics that do not pass through the vertex and meet 
each other transversely.
\end{enumerate}
\end{mtheorem}

\begin{proof}
(i) follows from \cite{Prokhorov2010a}. 

To prove (ii) we consider only the case \eqref{list-rho=2-degree=12}. All other cases are similar
(see \ref{remark-list-section-P3P3} for the case \eqref{list-section-P3P3}).
Then, in notation of \ref{definition-contractible-plane} and \ref{Corollary-nef-cone},
$\rho(\mathfrak X/ \mathfrak D)=2$ and the nef cone $\Nef(\mathfrak X/ \mathfrak D)$ 
has two edges. 
Thus there are
two extremal contractions 
$\mathfrak f_i: \mathfrak X\to \mathfrak Z_i$, $i=1$, $2$ over $\mathfrak D\ni 0$.
Let $\mathcal M_i$, $i=1$, $2$ be (integral) nef divisors generating edges of $\Nef(\mathfrak X/ \mathfrak D)$.
We can take $\mathcal M_i$ to be primitive elements of  $\Pic(\mathfrak X/ \mathfrak D)\simeq \ZZ^{\oplus 2}$.
Thus $\mathcal M_i=\mathfrak f_i^*\mathcal A_i$, where $\mathcal A_i$ is an   ample generator of
$\Pic(\mathfrak Z_i/ \mathfrak D)\simeq \ZZ$.
By Corollary \ref{Corollary-nef-cone} 
the map $\mathfrak f_i$ induces 
an extremal contraction $\mathfrak f_{i, s}: \mathfrak X_s\to \mathfrak Z_{i,s}$
on a each fiber. 
For a general fiber $\mathfrak X_s$, $s\neq 0$ we have
 $\mathfrak Z_{i,s}\simeq \PP^2$. 
Let $f_i: X\to Z_i$ be the contraction induced on the central fiber $X=\mathfrak X_0$
and let $M_i:=\mathcal M_i|_{X}$. Then 
$M_i=f_i^*A_i$, where $A_i$ is an ample divisor on $Z_i$.
By \eqref{rho2-eq-K-1} $M_1+M_2=-K_X$.
Since $M_i^3=(\mathcal M_i|_{\mathfrak X_s})^3=0$, $Z_i$ is a surface
and $f_i$ is a generically conic bundle.
Hence, $A_i^2=-\frac12 M_i^2\cdot K_X= -(\mathcal M_i|_{\mathfrak X_s})^2\cdot K_{\mathfrak X_s}=1$.
By semicontinuity $\dim H^0(Z_i, A_i)=\dim H^0(X, M_i)\ge \dim H^0(\mathfrak X_s, \mathcal M_i|_{\mathfrak X_s})=3$.
Hence, by \cite{Fujita-1975} $Z_i\simeq \PP^2$.
The map $f=f_1\times f_2: X\to \PP^2\times \PP^2$ must be finite  and $G$-equivariant.

Let $(d_1,d_2)$ be the  bidegree of $f(X)$ in $\PP^2\times \PP^2$ (as a divisor).
Then
\[
d_i=\OOO_{\PP^2\times \PP^2}(1,0)\cdot \OOO_{\PP^2\times \PP^2}(1,0)^2\cdot f(X) =
\frac{M_1\cdot M_2^2}{\deg f}=\frac{2}{\deg f}.
\]
Similarly, $d_2=2/\deg f$. Thus $d_1=d_2$ and there are two possibilities:
\begin{enumerate}
 \item[a)]
$f(X)$ is of bidegree $(2,2)$ and $f$ is birational;
 \item[b)]
$f(X)$ is of bidegree $(1,1)$ and $f$ is finite of degree $2$.
\end{enumerate}
In the first case, the map  $f:X\to f(X)$ 
is birational and finite.
By subadjunction $K_X=f^*K_{f(X)}-B$, where $B$ is an effective divisor.
On the other hand, $K_X=-M_1-M_2=f^*\OOO_{\PP^2\times \PP^2}(-1,-1)=f^*K_{f(X)}$.
Hence, $B=0$, $f(X)$ is normal, and $f:X\to f(X)$ is an isomorphism.

In the case b), we note that any irreducible singular hypersurface of 
bidegree $(1,1)$ in $\PP^2\times \PP^2$ 
can be given, in some coordinate system,  by the equation
$x_1y_1+x_2y_2=0$. In this case, the singular locus of $f(X)$ consists of one node and $f(X)$
contains exactly two planes. 
Therefore, $\rk \Cl(f(X))^G>1$ and so $\rk \Cl(X)>1$, a contradiction.
Hence $f(X)$ is smooth.

Now we prove (iii).
As above, by Corollary \ref{Corollary-nef-cone}
two extremal contractions 
$\mathfrak f_i: \mathfrak X\to \mathfrak Z_i$, $i=1$, $2$ over $\mathfrak D\ni 0$
induce
extremal contractions $\mathfrak f_{i, s}: \mathfrak X_s\to \mathfrak Z_{i,s}$
on a each fiber and for  $s\neq 0$ the variety  $\mathfrak Z_{i,s}$ is a  smooth quadric. 
Hence $\mathfrak f_i$ induces a birational extremal contraction 
$f_i: X\to Q_i$ on the central fiber $X=\mathfrak X_0$.
Let $\mathcal M_i=\mathfrak f_i^*\mathcal A_i$, where $\mathcal A_i$ is an   ample generator of
$\Pic(\mathfrak Z_i/ \mathfrak D)\simeq \ZZ$ and let $M_i:=\mathcal M_i|_{X}$. Then 
$M_i=f_i^*A_i$, where $A_i$ is an ample divisor on $Z_i$. 
By semicontinuity $\dim H^0(Q_i, A_i)\ge 5$.
Hence, by \cite{Fujita-1975} $Q_i$ is a 
quadric in $\PP^4$. In particular, $K_{Q_i}$ is Cartier.
Since $-K_X$ is ample, by standard facts of the minimal model program $Q_i$ has 
at worst terminal singularities.
Moreover, $f_i$ is an isomorphism over the singular locus of $Q_i$.
Hence $Q_i$ is either a smooth quadric or a quadratic cone in $\PP^4$.
Let $\mathcal E_i\subset \mathfrak X$ be $\mathfrak f_i$-exceptional divisor.
Since $\mathfrak X$ is factorial, $\mathcal E_i$ is an irreducible Cartier divisor.
Let $E_i:=\mathcal E_i|_X$. Clearly, the divisor $E_i$ coincides with the $f_i$-exceptional locus.
We have $-K_X=M_1+M_2=E_1+E_2$ (because this holds on a general fiber). 
 
Let $G_{\bullet}\subset G$ be the subgroup of index $2$ that 
stabilizes $M_1$ (and $M_2$). Then the contraction $f_i$ is $G_{\bullet}$-equivariant.
Since $\rk \Cl(Q_i)\le 2$, we have $\Cl(Q_i)^{G_{\bullet}}\simeq \ZZ$. 

Write $E_i=\sum_{j=1}^r k_iE_i^{(j)}$, where $k_i>0$, $r\ge 1$, and $E_i^{(j)}$ are prime Weil divisors.
Let $\Gamma_i:=f_i(E_i)$ and let $\Gamma_i^{(j)}=f_i(E_i^{(j)})$.
Since $\Cl(X)^G$ is generated by $K_X$, the divisors $E_i^{(j)}$ form one $G$-orbit and $k_i=1$ for all $i$.
In particular, $E_1$ and $E_2$ have no common components and 
$\Gamma_i^{(j)}$ is a curve for all $j$ (because 
$E_1\neq E_2$ and $E_1\cdot E_2\cdot M_i\neq 0$).
Let  $D\in |-K_X|$ be a general member and
let $D_i:=f_i(D)$. Then $D_i\in |-K_{Q_i}|$. 
Since the linear system $|-K_X|$ is base point free \cite{Jahnke-Radloff-2006},
$D$ is a smooth K3 surface. 
We may take $S$ so that it does not contain non-trivial fibers of $E_i\to \Gamma_i$
and meets a general fiber transversely at one point.
Then $D\to D_i$ is a finite birational morphism and $D_i$ has at worst isolated singularities.
Since $D_i$ is Cohen-Macaulay, it is smooth and $D\to D_i$
is an isomorphism. Then
\[
\begin{array}{ll}
\phantom{-}\deg \Gamma_i=A_i\cdot \Gamma_i=-K_X\cdot M_i\cdot E_i=4,  
\\[10pt]
-2=-K_X\cdot E_i^2=D\cdot E_i^2=\Gamma_i^2=2p_a(\Gamma_i)-2 \quad 
\Longrightarrow\quad p_a(\Gamma_i)=0.
\end{array}
\]
Since $(-K_X)^2\cdot E_1=\frac12 (-K_X)^3=14$, $r$ divides $14$.
Further, since $-K_X\cdot M_1\cdot  E_1=4$, $r$ divides $4$.

Therefore, $r=1$ or $2$. If $r=1$, i.e., $E_i$ is irreducible, then 
so $\Gamma_i:=f_i(E_i)$ is. Thus $\Gamma_i$ is a smooth rational curve
of degree $4$. In this case, 
$\Gamma_i$ is contained into the smooth locus of $Q_i$,
$f_i$ is the blowup of $\Gamma_i$, and
$X$ is smooth along $E_i$ \cite{Cutkosky-1988}. Assume that $Q_i$ is singular, i.e., 
it is the  projective cone over a smooth quadric $W\simeq \PP^1\times \PP^1$ in $\PP^3$.
Let $O_i\in Q_i$ be the vertex of the cone. Then $X$ has a unique singular point $O\in X$
and $f_i(O)=O_i$.  Moreover, $O\notin E_i$.
Since $E_1+E_2\sim -K_X$, $f_1(E_2)\sim -K_{Q_1}$.
Moreover, $f_1(E_2)$ is singular along $\Gamma_1$. Hence 
every  $2$-secant line to $\Gamma_1$  is contained in 
$f_1(E_2)$. Since $O_1\notin f_1(E_2)$, the projection of  $\Gamma_1$
from $O_1$ to $W$ is an isomorphism. Thus $\Gamma_1$ is contained into a divisor of type 
$(1,3)$ on $Q_1$.
But in this case, the action of $G_{\bullet}$ on  $\Cl(Q_i)\simeq \ZZ\oplus\ZZ$
must be trivial, a contradiction.

Finally, assume that $r=2$, i.e., $E_i=E_i^{(1)}+E_i^{(2)}$. 
Since the divisors $E_i^{(1)}$ and $E_i^{(2)}$ are not Cartier, 
$G_\bullet$  interchanges $E_i^{(1)}$ and $E_i^{(2)}$.
Hence, $G_\bullet$  interchanges also $f_1(E_2^{(1)})$ and $f_1(E_2^{(2)})$.
On the other hand, $f_1(E_2^{(1)})+f_1(E_2^{(2)})\sim -K_{Q_1}\sim \OOO_{Q_1}(3)$.
This is possible only if $f_1(E_2^{(1)})$ and $f_1(E_2^{(2)})$
are not Cartier. So, $Q_1$ is the quadratic cone and $f_1(E_2)$ contains its vertex.
Again by \cite{Cutkosky-1988} $X$ is the blowup of $Q_1$ along $\Gamma_1$.
\end{proof}


\def\cprime{$'$} \def\mathbb#1{\mathbf#1}
  \def\bblapr{April}\def\mathbb#1{\mathbf#1}

\end{document}